\documentclass[12pt, reqno]{amsart}
\usepackage{amsmath}
\usepackage{amssymb}

\numberwithin{equation}{section}

\begin{document}
\setlength{\baselineskip}{1.30em}

{\theoremstyle{plain}
	\newtheorem{Theorem}{\bf Theorem}[section]
	\newtheorem{Proposition}[Theorem]{\bf Proposition}
	\newtheorem{claim}[Theorem]{\bf Claim}
	\newtheorem{Lemma}[Theorem]{\bf Lemma}
	\newtheorem{Corollary}[Theorem]{\bf Corollary}
}
{\theoremstyle{remark}
	\newtheorem{Remark}[Theorem]{\bf Remark}
	\newtheorem{Example}[Theorem]{\bf Example}
}
{\theoremstyle{definition}
	\newtheorem{defn}[Theorem]{\bf Definition}
}

%%%%%%%%%%%%%%%%%%%%%%%%%%%%%%%%%%%%%%%%%%%

\def\To{\longrightarrow}
\def\height{\operatorname{ht}}
\def\reg{\operatorname{reg}}
\def\Hom{\operatorname{Hom}}
\def\Proj{\operatorname{Proj}}
\def\BiProj{\operatorname{BiProj}}
\def\grade{\operatorname{grade}}
\def\Spec{\operatorname{Spec}}

\def\mm{{\frak m}}
\def\pp{{\frak p}}
\def\O{{\mathcal O}}
\def\I{{\mathcal I}}
\def\M{{\mathcal M}}
\def\L{{\mathcal L}}
\def\S{{\mathcal S}}
\def\NN{{\mathbb N}}
\def\PP{{\mathbb P}}
\def\ZZ{{\mathbb Z}}
\def\T{{\mathcal T}}

%%%%%%%%%%%%%%%%%%%%%%%%%%%%%%%%%%%%%%%%%%%

\title{Asymptotic behaviour of arithmetically Cohen-Macaulay blow-ups}
\author{Huy T\`ai H\`a}
\author{Ng\^o Vi\^et Trung}
\thanks{The second author is partially supported by the National Basic 
Research Program of Vietnam} 
\address{Department of Mathematics, University of Missouri, Columbia MO 65201, USA}
\email{tai@math.missouri.edu}
\address{Institute of Mathematics, 18 Hoang Quoc Viet, Hanoi, Vietnam}
\email{nvtrung@math.ac.vn}

\begin{abstract}
This paper addresses problems related to the existence of arithmetic Macaulayfications of projective schemes. Let $Y$ be the blow-up of a projective scheme $X = \Proj R$ along the ideal sheaf of $I \subset R$. It is known that there are embeddings $Y \cong \Proj k[(I^e)_c]$ for $c \ge d(I)e + 1$, where $d(I)$ denotes the maximal generating degree of $I$, and that there exists 
a Cohen-Macaulay ring of the form $k[(I^e)_c]$ if and only if $H^0(Y,\O_Y) = k$, $H^i(Y,\O_Y) = 0$ for $i = 1,...,\dim Y-1$, $Y$ is equidimensional and Cohen-Macaulay. Cutkosky and Herzog asked when there is a linear bound on $c$ and $e$ 
ensuring that $k[(I^e)_c]$ is a Cohen-Macaulay ring. We obtain a surprising compelte answer to this question, namely, that under the above conditions, there are well determined invariants $\varepsilon$ and $e_0$ such that $k[(I^e)_c]$  is Cohen-Macaulay for all $c > d(I)e + \varepsilon$ and $e > e_0$.  Our approach is based on recent results on the asymptotic linearity of the Castelnuovo-Mumford regularity of ideal powers. We also  investigate the existence of  a Cohen-Macaulay Rees algebra of the form $R[(I^e)_ct]$ (which provides an arithmetic Macaulayfication for $X$). If $R$ has negative $a^*$-invariant, we prove that such a Cohen-Macaulay Rees algebra exists if and only if $\pi_*\O_Y = \O_X$, $R^i\pi_*\O_Y = 0$ for $i > 0$, $Y$ is equidimensional and Cohen-Macaulay. Especially, these conditions imply the Cohen-Macaulayness of  $R[(I^e)_ct]$  for all $c > d(I)e + \varepsilon$ and $e > e_0$.  The above results can be applied to obtain several new classes of Cohen-Macaulay algebras. \end{abstract}
\maketitle

%%%%%%%%%%%%%%%%%%%%%%%%%%%%%%%%%%%%%%%%%%%

\section*{Introduction} \smallskip

\noindent Let $X$ be a projective scheme over a field $k$.
An arithmetic Macaulayfication of $X$ is a proper birational morphism $\pi: Y \to X$ such that $Y$ has an arithmetically Cohen-Macaulay embedding, i.e. there exists a Cohen-Macaulay standard graded algebra $A$ over $k$ such that $Y \cong \Proj A$.
Inspired by the problem of desingularization, one may ask when $X$ has an arithmetic Macaulayfication. This problem is a global version of the problem of arithmetic Macaulayfication of local rings recently solved by Kawasaki \cite{ka}. The existence of an arithmetic Macaulayfication is usually obtained by blowing up $X$ at a suitable subscheme. \par

Let $R$ be a standard graded $k$-algebra and $I \subset R$ a homogeneous ideal such that $X = \Proj R$ and $Y$ is the blow-up of $X$ with respect to the ideal sheaf of $I$. It was observed by Cutkosky and Herzog \cite{che} that $Y \cong \Proj k[(I^e)_c]$  for $c \ge d(I)e + 1$,  where $(I^e)_c$ denotes the vector space of forms of degree $c$ of the ideal power $I^e$ and $d(I)$ is the maximal degree of the elements of a homogeneous basis of $I$. In other words, $Y$ can be embedded into a projective space by the complete linear system $|cE_0 - eE|$, where $E$ denotes the exceptional divisor and $E_0$ is the pull-back of a general hyperplane in $X$.  By \cite{v2} we know that there exists a Cohen-Macaulay ring $k[(I^e)_c]$ for $c \ge d(I)e + 1$ if and only if $Y$ satisfies the following conditions:\par
\begin{itemize}
\item  $Y$ is equidimensional and Cohen-Macaulay, 
\item  $H^0(Y,\O_Y) = k$ and $H^i(Y,\O_Y) = 0$ for $i = 1,...,\dim Y-1$.
\end{itemize}

In the first part of this paper, we study the problem for which values of $c$ and $e$ is $k[(I^e)_c]$ a Cohen-Macaulay ring. This problem originated from  a beautiful result of Geramita, Gimigliano and Pitteloud \cite{ggp} which shows that if $I$ is the defining ideal of a set of fat points in a projective space over a field of characteristic zero, then 
$k[I_c]$ is a Cohen-Macaulay ring for all $c \ge \reg(I)$, where $\reg(I)$ is the Castelnuovo-Mumford regularity of $I$. This result initiated the study on the Cohen-Macaulayness of algebras of the form $k[(I^e)_c]$ first in \cite{chtv} and then in \cite{che, v1, v2, ha2}. In particular, Cutkosky and Herzog \cite{che} showed that if $I$ is a locally complete intersection ideal, then there exists a constant $\delta$ such that $k[(I^e)_c]$ is Cohen-Macaulay for $c \ge \delta e$. They asked when there is a linear bound on $c$ and $e$ ensuring that $k[(I^e)_c]$ is a Cohen-Macaulay ring. 

Our results will give a complete answer to this question.
We show that if the above two conditions are satisfied, then there exist well-determined invariants $\varepsilon$ and $e_0$  such that 
$k[(I^e)_c]$ is a  Cohen-Macaulay ring for all $c > d(I)e + \varepsilon$ and $e > e_0$ (Theorem \ref{diagonal-1}). The invariant $e_0$ is a projective version of the $a^*$-invariant, which is the largest non-vanishing degree of the graded local cohomology modules \cite{sh, t1}. The invariant $\varepsilon$ comes from the asymptotic linearity of the Castelnuovo-Mumford regularity of powers of ideals (\cite{sw, cht,  k, tw}). We will see that the bounds $c > d(I)e + \varepsilon$ and $e > e_0$ are the best possible for the existence of a Cohen-Macaulay ring $k[(I^e)_c]$ (Proposition \ref{sharp-1} and Example \ref{sharp-example}). 
In particular, if the Rees algebra $R[It]$ is locally Cohen-Macaulay on $X$, then $e_0 = 0$ and we may replace the second condition by the weaker condition that $H^0(X,\O_X) = k$ and $H^i(X,\O_X) = 0$ for $i = 1,...,\dim X-1$ (Theorem \ref{diagonal-2}).   These results unify all previously known results on the Cohen-Macaulayness of $k[(I^e)_c]$ which were obtained by different methods.   \par

In the second part of this paper, we investigate the more difficult question of when $Y$ is an arithmetically Cohen-Macaulay blow-up of $X$; that is, when there exists a standard graded $k$-algebra $R$ and an ideal $J \subset R$, such that $X = \Proj R$, $Y$ is the blow-up of $X$ along the ideal sheaf of $J$, and $R[Jt]$ is a Cohen-Macaulay ring.  Given $R$ and $I$, we will concentrate on ideals $J  \subseteq I$  which are generated by the elements of $(I^e)_c$. It is obvious that $I^e$ and $J$ define the same ideal sheaf for $c \ge d(I)e$. Rees algebras of the form  $R[I_ct]$ ($e = 1$) have been studied first for the defining ideal of a set of points in \cite{ha1} and then for locally complete intersection ideals in \cite{cha}, where it was shown that there exists a constant $\lambda$ such that $R[I_ct]$ is a Cohen-Macaulay ring for $c \ge \lambda$. This leads to the problem of whether there is a constant $\delta$ such that the Rees algebra $R[(I^e)_ct]$ is a Cohen-Macaulay ring for $c \ge \delta e$. 

If $a^*(R) < 0$ (e.g. if $R$ is a polynomial ring) we solve this problem by showing that there exists a Cohen-Macaulay ring $R[(I^e)_ct]$ with $c \ge d(I)e$ if and only if the following conditions are satisfied:
\begin{itemize}
\item $Y$ is equidimesional and Cohen-Macaulay,
\item $\pi_*\O_Y = \O_X$, $R^i\pi_*\O_Y = 0$ for $i > 0$.
\end{itemize}
Especially, these conditions imply that $R[(I^e)_ct]$ is a Cohen-Macaulay ring for all $c > d(I)e+ \varepsilon$ and $e > e_0$ (Theorem \ref{truncated-1}). From this it follows that there exists a Cohen-Macaulay algebra of the form $R[I_ct]$ with $c \ge d(I)$ if and only if $R[It]$ is locally Cohen-Macaulay on $X$ and that $e_0 = 0$ in this case (Corollary \ref{e=1}). We would like to point out that this phenomenon does not hold in general. In fact, there exist  examples with  $a^*(R) \ge 0$ such that $R[(I^e)_{d(I)e}t]$ is a Cohen-Macaulay ring, whereas $R[(I^e)_ct]$ is not a Cohen-Macaulay ring for any $c > d(I)e$ (Example \ref{non-linear}). 
Using the above result we obtain several new classes of Cohen-Macaulay Rees algebras. Furthermore, we show that if $H^0(X,\O_X) = k$ and $H^i(X,\O_X) = 0$ for $i > 0$, then $Y$ is an arithmetically Cohen-Macaulay blow-up of $X$ if and only if $Y$ is locally arithmetical Cohen-Macaulay on $X$ (Theorem \ref{truncated-2}). \par

Our approach is based on the facts that  the Rees algebra $S = R[It]$ has a natural bi-gradation and that $k[(I^e)_c]$ can be viewed as a diagonal subalgebra of $S$ \cite{chtv}. As a consequence, the Cohen-Macaulayness of $k[(I^e)_c]$ can be characterized by means of the sheaf cohomology $H^i(Y,\O_Y(m,n))$. Using Leray spectral sequence and Serre-Grothendieck correspondence, we may pass this sheaf cohomology to the local cohomology of $I^n$ and of $\omega_n$, where $\omega_S = \oplus_{n \in \ZZ}\omega_n$ denotes the graded canonical module of $S$. 
It was shown recently that there are linear bounds for the vanishing of the local cohomology of $I^n$ and $\omega_n$ (\cite{sw, cht, k, tw}). It turns out that these linear bounds yield a linear bound on $c$ and $e$ such that $k[(I^e)_c]$ is a Cohen-Macaulay ring. The Cohen-Macaulayness of  the Rees algebra $R[(I^e)_ct]$ can be studied similarly by using a recent result of Hyry \cite{hy} which characterizes the Cohen-Macaulayness of a standard bi-graded algebra by means of sheaf cohomology.\par

The paper is organized as follows. In Section 1, we introduce the notion of a projective $a^*$-invariant which governs how sheaf cohomology behaves through blow-ups.  In Section 2, we study the Cohen-Macaulayness of rings of the form $k[(I^e)_c]$ which correspond to projective embeddings of $Y$. The last section of the paper deals with the problem of when $Y$ is an arithmetically Cohen-Macaulay blow-up of $X$. \par

For unexplained notations and facts we refer the reader to the books \cite{bs, bh, har}.

\section{$a^*$-invariants}

\noindent Let $R$ be an arbitrary commutative noetherian ring. Let 
$S = \oplus_{n \ge 0}S_n$ be a finitely generated graded algebra over $R$. We shall always use $S_+ = \oplus_{n > 0}S_n$ to denote the ideal generated by the homogeneous elements of positive degrees of $S$. Given any finitely generated graded $S$-module $F$, the local cohomology module $H_{S_+}^i(F)$ is also a graded $S$-module. It is well-known that $H_{S_+}^i(F)_n = 0$ for  $n \gg 0$, $i \ge 0$. Put
$$a_i(F) = \left\{ \begin{array}{lll}  -\infty & \text{if} & H_{S_+}^i(F) = 0,\\
\max\{n|\ H_{S_+}^i(F)_n \neq 0\} & \text{if} &H_{S_+}^i(F) \neq 0. \end{array} \right. $$
Note that $a(F) := a_{\dim F}(F)$ is called the $a$-invariant of $F$ if $S$ is a standard graded algebra over a field. The {\it $a^*$-invariant} of $F$ is defined to be
$$a^*(F) := \max\{a_i(F)|\ i \ge 0\}.$$
This invariant was introduced in \cite{t1} and \cite{sh} in order to control the vanishing of graded local cohomology modules with different supports. 
It is closely related to the Castelnuovo-Mumford regularity via the equality
$$\reg(F) = \max\{a_i(F)+i|\ i \ge 0\}.$$ \par

Here we are interested in the case when $R$ is a standard graded algebra over a field $k$ and $S = R[It]$ is the Rees algebra of a homogeneous ideal $I \subset R$ with $\height I \ge 1$. This Rees algebra has a natural grading with $S_n = I^nt^n$. Let $\omega_S = \oplus_{n \in \ZZ}\omega_n$ denote the canonical graded module of $S$.

\begin{Lemma} \label{Rees}
Let $S = R[It]$ be as above. If $S$ is a Cohen-Macaulay ring, then $a^*(S) = -1$ and $a^*(\omega_S)  = 0$.
\end{Lemma}

\begin{proof}
It is well-known that $\dim S = \dim R+1$. Since $S/S_+ = R$, we have $\height S_+ = \dim S - \dim R = 1$. This implies $\grade S_+ = 1$.
Hence $a^*(S) \ge -1$ by \cite[Corollary 2.3]{t1}. On the other hand, the Cohen-Macaulayness of $S$ implies 
$H_M^i(S) = 0$ for $i < \dim S$, where $M$ denotes the maximal graded ideal of $S$. By  \cite[Corollary 3.2]{ti} we always have $H_M^{\dim S}(S)_n = 0$ for $n \ge 0$. Hence $H_M^i(S)_n = 0$ for all $n \ge 0$ and $i  \ge 0$. By \cite[Lemma 2.3]{hy} (or \cite[Corollary 2.8]{t1}), this implies $H_{S_+}^{i}(S)_n = 0$ for all $n \ge 0$ and $i \ge 0$. Therefore, $a^*(S) = -1$.  \par
Since $\omega_S$ is a Cohen-Macaulay module with $\Hom_S(\omega_S,\omega_S) \cong S$  \cite[Proposition 2]{ao}, we also have $H_M^i(\omega_S) = 0$ for $i < \dim S$ and, by local duality, 
$$H_M^{\dim S}(\omega_S)_n \cong \Hom_S(\omega_S,\omega_S)_{-n} \cong S_{-n}.$$
Since $S_0 = R \neq 0$ and $S_{-n} = 0$ for $n > 0$, we can conclude that $a_X^*(\omega_S) = 0$.
\end{proof}

Let $X = \Proj R$. For each $\pp \in X$, the homogeneous localization $F_{(\pp)}$ is a finitely generated graded module over $S_{(\pp)}$. Hence, we can define the {\it projective $a^*$-invariant}
$$a_X^*(F) := \max\{a^*(F_{(\pp)})|\ \pp \in X\}.$$
Note that $H^i_{S_{(\pp)+}}(F_{(\pp)}) = H_{S_+}^i(F)_{(\pp)}$ (cf. \cite[Remark 2.2]{sh}). Then we always have $a_X^*(F) \le a^*(F).$ Hence $a_X^*(F)$ is  a finite number. Since $a_X^*(F)$ is determined by the local structure of $F$ on $X$, it can easily be estimated in certain situations. As a demonstration, we show how to estimate $a_X^*(F)$ in the following case which will play an important role in our further investigation. \par

We say that $S$ is {\it locally Cohen-Macaulay on $X$} if $S_{(\pp)} $ is a Cohen-Macaulay ring for every $\pp \in \Proj R$. This condition holds if, for instance, $X$ is locally Cohen-Macaulay and $\I$ is locally a complete intersection.\par

\begin{Proposition} \label{arithmetic}
Let $X = \Proj R$ and $S = R[It]$ be as above. Then $a_X^*(S) \ge -1$ and $a_X^*(\omega_S) \ge 0$. Equalities hold if $S$ is locally Cohen-Macaulay on $X$. 
\end{Proposition}

\begin{proof}
Let $\pp$ be a minimal prime ideal in $X$. Then $R_{(\pp)}$ is an artinian ring. Since $\pp \not\supseteq I$, we have $I_{(\pp)} = R_{(\pp)}$. Hence $S_{(\pp)} = R_{(\pp)}[t]$ is a Cohen-Macaulay ring.  By Lemma \ref{Rees}, this implies 
$a^*(S_{(\pp)}) = -1$ and $a^*(\omega_{S_{(\pp)}}) = 0$. Hence
$a_X^*(S) \ge -1$ and $a_X^*(\omega_S) \ge 0$. This proves the first statement. The second statement is an immediate consequence of Lemma \ref{Rees}.
\end{proof}

Beside the natural $\NN$-graded structure given by the degrees of $t$, the Rees algebra $S = R[It]$ also has a natural bi-gradation with 
$$S_{(m,n)} = (I^n)_mt^n$$
for $(m,n) \in \NN^2$. 
Let $Y$ be the blow-up  of $X$ along the ideal sheaf of $I$. Then  $Y = \Proj S$ with respect to this bi-gradation. If $F = \oplus_{(m,n) \in \ZZ^2}F_{(m,n)}$ is a finitely generated bi-graded $S$-module, then $F$ is also an $\ZZ$-graded $S$-module with $F_n = \oplus_{m \in \ZZ}F_{(m,n)}$. Let $\widetilde F$ denote the sheaf associated to $F$ on $Y$. We write $\tilde{F}(n)$ and $\tilde{F}(m,n)$ to denote the twisted $\O_Y$-modules with respect to the $\NN$-gradation and the $\NN^2$-gradation of $S$. Moreover, we denote by $\widetilde{F_n}$ the sheafification of $F_n$ on $X$.  \par

It turns out that $a^*_X(F)$ is a measure for when we can pass from the sheaf cohomology of $\widetilde F(m,n)$ on $Y$ to that of $\widetilde{F_n}(m)$ on $X$.

\begin{Proposition} \label{Leray}
Let $F$ be a finitely generated bi-graded $S$-module. For $n > a_X^*(F)$ we have\par
{\rm (i)} $\pi_*(\widetilde F(n)) = \widetilde{F_n}$ and $R^i\pi_*(\widetilde F(n)) = 0$ for $i > 0$,\par
{\rm (ii)} $H^i(Y,\widetilde F(m,n)) \cong H^i(X,\widetilde{F_n}(m))$ for all $m \in \ZZ$ and $i\ge 0$.
\end{Proposition}

\begin{proof} 
Since (i) is a local statement, we only need to show that it holds locally. 
Let $\pp$ be a closed point of $X$, and consider the restriction $\pi_\pp$ of $\pi$ over an affine open neighborhood $\Spec \O_{X, \pp}$ of $\pp$
$$ \pi_\pp : Y_\pp = Y \times_X \Spec \O_{X, \pp} \to \Spec \O_{X, \pp}. $$
We have $\tilde{F} \big|_{Y_\pp} = \widetilde{F_{(\pp)}}$, where $\widetilde{F_{(\pp)}}$ is the sheaf associated to $F_{(\pp)}$ on $Y_\pp$. Thus,
$$R^i \pi_* (\tilde{F}(n)) \Big|_{\Spec \O_{X, \pp}} = R^i {\pi_\pp}_* (\widetilde{F_{(\pp)}}(n)) = H^i(Y_\pp, \widetilde{F_{(\pp)}}(n))\widetilde{~}. $$
On the other hand, we know by the Serre-Grothendieck correspondence that there are the exact sequence
$$ 0 \to H^0_{S_{(\pp)+}}(F_{(\pp)})_n \to (F_{(\pp)})_n \to H^0(Y_\pp, \widetilde{F_{(\pp)}}(n)) \to H^1_{S_{(\pp)+}}(F_{(\pp)})_n \to 0 $$
and the isomorphisms $H^i(Y_\pp, \widetilde{F_{(\pp)}}(n)) \cong H^{i+1}_{S_{(\pp)+}}(F_{(\pp)})_n$ for $i > 0$. By the definition of $a_X^*(F)$, we know that $H^i_{S_{(\pp)+}}(F_{(\pp)})_n$ for $n > a_X^*(F)$, $i > 0$. Thus,
$$R^i \pi_* (\tilde{F}(n)) \Big|_{\Spec \O_{X, \pp}} = H^i(Y_\pp,\widetilde{F_{(\pp)}}(n))\widetilde{~} = 
\left\{ \begin{array}{lll}  \widetilde{(F_n)_{(\pp)}} & \text{if} & i = 0,\\
0 & \text{if} &  i > 0, \end{array}\right. $$
for $n > a^*_X(F)$.   \par

To show (ii) we first observe that $\widetilde F(m,n) =  \widetilde F(n) \otimes  \pi^*\O_X(m)$.  By the projection formula, we have
$$R^i\pi_*(\widetilde F(m,n)) = R^i\pi_*(\widetilde F(n)) \otimes \O_X(m)
= \left\{ \begin{array}{lll} 
\widetilde{F_n}(m) & \text{if} &  i = 0,\\
0 & \text{if} & i > 0, \end{array} \right. $$
Hence the conclusion follows from the Leray spectral sequence 
$$H^i(X, R^j \pi_* (\tilde{F}(m,n))) \Rightarrow H^{i+j}(Y, \tilde{F}(m,n)).$$ 
\end{proof}

Let $Y$ be the blow-up of a projective scheme $X$ along an ideal sheaf $\I$. We say that $Y$ is {\it locally arithmetic Cohen-Macaulay on $X$} if there exist $R$ and $I$ such that $X = \Proj R$, $\I = \widetilde{I}$ and $S = R[It]$ is locally Cohen-Macaulay on $X$. 

\begin{Corollary} \label{transfer}
Assume that $Y$ is locally arithmetic Cohen-Macaulay on $X$. Then\par
{\rm (i) } $\pi_*\O_Y = \O_X$ and $R^i\pi_*\O_Y = 0$ for $i > 0$,\par
{\rm (ii) } $H^i(Y,\O_Y(m,0)) \cong H^i(X, \O_X(m))$ for  all $m \in \ZZ$, $i \ge 0$. 
\end{Corollary}

\begin{proof}
With the above notations we have $a_X^*(S) = -1$ by Lemma \ref{arithmetic}. Hence the conclusion follows from Proposition \ref{Leray} by putting $F = S$ and $n = 0$.
\end{proof}

For each $n$, the graded $R$-module $F_n$ has an $a^*$-invariant $a^*(F_n)$, which controls the vanishing of $H^i(X,\widetilde{F_n}(m))$ by the Grothendieck-Serre correspondence. On the other hand,  since $F$ is a finitely generated graded module over $S = R[It]$, there exists a number $n_0$ such that $F_n = I^{n-n_0}F_{n_0}$ for $n \ge n_0$. It  was recently discovered that for any finitely generated graded $R$-module $E$, the Castelnuovo-Mumford regularity 
$\reg (I^nE) $
is bounded by a linear  function on $n$ with slope $d(I)$ \cite[Theorem 2.2]{tw} (see also \cite{cht, k} for the case $R$ is a polynomial ring). By definition, we always have $$a^*(I^nE) \le \max\{a_i(I^nE)+i|\ i \ge 0\} = \reg(I^nE).$$
Therefore, $a^*(F_n)$ is bounded above by a linear function of the form $d(I)n + \varepsilon$ for $n \ge 1$. \par

We will denote by $\varepsilon(I)$ the smallest non-negative number such that 
$$a^*(I^n) \le d(I)n+\varepsilon(I)$$
for all $n \ge 1$. Since $\omega_S = \oplus_{n\in \ZZ}\omega_n$ is a finitely generated bi-graded $S$-module, there is a similar bound for $a^*(\omega_n)$. 
Note that the $R$-graded module $\omega_n$ is also called an {\it adjoint-type module} of $I$ because of its relationship to the adjoint ideals \cite{hs}. We will denote by $\varepsilon^*(I)$ the smallest non-negative number such that 
$$a_i(\omega_n) \le d(I)n+\varepsilon^*(I)$$
for $i \ge 2$ and $n \ge 1$. \par

The meaning of these invariants will become more apparent in the next sections.
Here we content ourselves with the following observations.

\begin{Lemma} \label{asymptotic}
With the above notations we have\par 
{\rm (i) } $H^0(X,\widetilde{S_n}(m))  = S_{(m,n)}$ and
$H^i(X,\widetilde{S_n}(m)) = 0$ for $ i > 0$ and $m  > d(I)n+\varepsilon(I)$,\par
{\rm (ii) } $H^i(X,\widetilde{\omega_n}(m)) = 0$ for $ i > 0$ and $m  > d(I)n+\varepsilon^*(I)$. 
\end{Lemma}

\begin{proof}
Since $S_n \cong  I^n$, we have $H_{R_+}^i(S_n)_m = 0$ for $i \ge 0$, $m  > d(I)n+\varepsilon(I)$ and $n \ge 1$. Hence the first statement follows from  the Serre-Grothendieck correspondence, which gives the exact sequence
$$0 \to H_{R_+}^0(S_n)_m \to S_{(m,n)} \to H^0(X,\widetilde{S_n}(m)) \to H_{R_+}^1(S_n)_m \to 0$$
and the isomorphisms 
$$H^i(X,\widetilde{S_n}(m)) \cong H_{R_+}^{i+1}(S_n)_m$$
for $i > 0$. The second statement can be proved similarly.
\end{proof}

\section{Arithmetically Cohen-Macaulay embeddings of blow-ups}

Let $X$ be a projective scheme over a field $k$. Let $Y \to X$ be the blowing up of $X$ along an ideal sheaf $\I$. We say that $Y$ has an {\it arithmetically Cohen-Macaulay embedding} if there exists a Cohen-Macaulay standard graded $k$-algebra $A$ such that $Y \cong \Proj A$.

Let $R$ be a finitely generated standard graded $k$-algebra and $I \subset R$ a homogeneous ideal such that $X = \Proj R$ and $\I$ is the ideal sheaf associated to $I$.  Let $S = R[It]$ be the Rees algebra of $R$ with respect to $I$. It is well-known that  $Y \cong \Proj k[(I^e)_c]$ for $c \ge d(I)e+1$ and $e \ge 1$, where $k[(I^e)_c]$ is the algebra generated by all forms of degree $c$ of the ideal power $I^e$ and $d(I)$ denotes the largest degree of a minimal set of homogeneous generators of $I$. (cf. \cite[Lemma 1.1]{che}). There is the following simple criterion for the existence of a Cohen-Macaulay algebra $k[(I^e)_c]$ (which is at the same time a criterion for the existence of an arithmetically Cohen-Macaulay embedding).

\begin{Lemma} \label{criterion} \cite[Corollary 3.5]{v2} 
There exists a Cohen-Macaulay ring $k[(I^e)_c]$ for $c \ge d(I)e+1$ if and only if the following conditions are satisfied: \par
{\rm (i) } $Y$ is equidimensional and Cohen-Macaulay,\par
{\rm (ii) } $H^0(Y,\O_Y) = k$ and $H^i(Y,\O_Y) = 0$ for $i = 1,...,\dim Y-1$.
\end{Lemma}

The proof of \cite{v2} used a deep result on the relationship between the local cohomology modules of a bi-graded algebra and its diagonal subalgebras \cite{chtv}. 
However, the above lemma simply follows from the basic fact that (i) and (ii) are equivalent to the existence of an arithmetically Cohen-Macaulay Veronese embedding of $Y$, (cf. \cite[Lemma 1.1]{cha}). In fact, the Veronese subalgebras of $k[I_c]$ are exactly the algebras of the form $k[(I^e)_{ce}]$ for $c \ge d(I)+1$, $e \ge 1$. We notice that the statements of \cite[Corollary 3.5]{v2} and \cite[Lemma 1.1]{cha} missed the equidimensional condition. 

In this section we will determine for which values of $c$ and $e$ is $k[(I^e)_c]$ a Cohen-Macaulay ring. First, we show that there are well determined invariants $\varepsilon$ and $e_0$ such that $k[(I^e)_c]$ is a Cohen-Macaulay ring for all $c > d(I)e+ \varepsilon$ and $ e > e_0$. \par

\begin{Theorem} \label{diagonal-1}
Let $R$ be a standard graded algebra over a field $k$ and $I \subset R$ a homogeneous ideal with $\height I \ge 1$.  Let $Y$ be the blow-up of $X = \Proj R$ along the ideal sheaf of $I$ and $S = R[It]$. Assume that \par
{\rm (i) } $Y$ is equidimensional and Cohen-Macaulay,\par
{\rm (ii) } $H^0(Y,\O_Y) = k$ and $H^i(Y,\O_Y) = 0$ for $i = 1,...,\dim Y-1$.\par
\noindent Then $k[(I^e)_c]$ is a Cohen-Macaulay ring for  $c > d(I)e + \max\{\varepsilon(I),\varepsilon^*(I)\}$ and $e > \max\{a_X^*( S),a_X^*(\omega_S)\}$. 
\end{Theorem}

Note first  that we always have 
$\max \{ a^*_X(S), a^*_X(\omega_S) \} \ge 0$  by Proposition \ref{arithmetic} and $\max\{\varepsilon(I),\varepsilon^*(I)\} \ge 0$ by the definition of $\varepsilon(I)$ and $\varepsilon^*(I)$.

\begin{proof}
Let $A = k[(I^e)_c]$.  Since $c \ge de + 1$, we have $Y \cong \Proj A$ \cite[Lemma 1.1]{che}. On the other hand,  the Rees algebra $S = R[It]$ has a natural bi-gradation with $S_{(m,n)} = (I^n)_mt^n$ and $Y = \Proj S$.  Moreover, we may view $A$ as a diagonal subalgebra of $S$; that is,
$A = \oplus_{n \in \NN}S_{(cn,en)}$ \cite[Lemma 1.2]{chtv}. From this it follows that $A(n) \tilde{~} = \O_Y(cn,en)$. Therefore, the Serre-Grothendieck correspondence yields the exact sequence
$$0 \To H_{A_+}^0(A) \To A \To \oplus_{n \in \ZZ}H^0(Y,\O_Y(cn,en)) \To H_{A_+}^1(A) \To 0$$
and the isomorphisms
$$\oplus_{n \in \ZZ}H^i(Y,\O_Y(cn,en)) \cong H_{A_+}^{i+1}(A)$$
for $i \ge 1$. 
It is well-known that $A$ is a Cohen-Macaulay ring if and only if $H_{A_+}^i(A) = 0$ for $i \neq \dim A$.  Therefore, $A$ is a Cohen-Macaulay ring if we can show \begin{align*}
H^0(Y,\O_Y(cn,en)) & = A_n = 
\left\{\begin{array}{lll} 0 & \text{for} & n < 0,\\ k & \text{for} & n = 0,\\
(I^{en})_{cn} & \text{for} & n > 0,\end{array}\right.\\
H^i(Y,\O_Y(cn,en)) & = 0\ (i = 1,...,\dim Y-1).
\end{align*}

For $n = 0$, this follows from the assumption $H^0(Y,\O_Y) = k$ and
$H^i(Y,\O_Y) = 0$ for $i = 1,...,\dim Y-1$. \par

For $n > 0$ we have $cn > d(I)en+\varepsilon(I)n \ge d(I)en+\varepsilon(I)$ and $en > a_X^*(S)n \ge a_X^*(S)$. Hence, using Proposition \ref{Leray} and Lemma \ref{asymptotic} we get
\begin{align*}
H^0(Y,\O_Y(cn,en)) & = H^0(X,\widetilde{I^{en}}(cn)) = (I^{en})_{cn},\\
H^i(Y,\O_Y(cn,en)) & = H^i(X,\widetilde{I^{en}}(cn)) = 0,\  i = 1,...,\dim Y-1.
\end{align*}\par

For $n < 0$ we have 
$$
H^i(Y,\O_Y(cn,en)) =  H^{\dim Y-i}(Y,\omega_Y(-cn,-en))
$$
for $i \ge 0$. Serre duality  can be applied here because $Y$ is equidimensional and Cohen-Macaulay. Since $-cn > -d(I)en - \varepsilon^*(I)n \ge -d(I)en + \varepsilon^*(I)$ and $-en > - a_X^*(\omega_S)n \ge a_X^*(\omega_S)$, using Proposition \ref{Leray} and Lemma \ref{asymptotic} we get
$$H^{\dim Y-i}(Y,\omega_Y(-cn,-en)) = H^{\dim Y-i}(X,\widetilde{(\omega_S)_{-en}}(-cn)) = 0$$
for $i < \dim Y$. So we get $H^i(Y,\O_Y(cn,en)) = 0$ for all $n < 0$ and $i = 0,...,\dim Y-1$. The proof of Theorem \ref{diagonal-1} is now complete.
\end{proof}

The following proposition shows that the bound $e > \max \{ a^*_X(S), a^*_X(\omega_S) \}$ of Theorem \ref{diagonal-1} is the best possible.

\begin{Proposition} \label{sharp-1}
Let  the notations and assumptions be as in Theorem \ref{diagonal-1}.  Put 
$$e_0 = \max \{a^*_X(S), a^*_X(\omega_S)\}.$$
Then $k[(I^{e_0})_c]$ is not a Cohen-Macaulay ring for $c \gg 0$ if $e_0 \ge 1$.
\end{Proposition}

\begin{proof} 
Let $A = k[(I^{e_0})_c]$ for $c \gg 0$. As we have seen in the proof of Theorem \ref{diagonal-1}, $A$ is not Cohen-Macaulay if $H^0(Y, \O_Y(c,e_0)) \not= (I^{e_0})_c$ or $H^i(Y, \O_Y(c,e_0)) \not= 0$ or $H^i(Y,\O_Y(-c,-e_0)) \neq 0$ for some $i = 1,...,\dim Y-1$. \par

We shall first consider the case $e_0 = a^*_X(S)$. 
Let $q$ be the smallest integer such that 
$e_0 = \max\{a_q(S_{(\pp)})|\ \pp \in X\}.$ Then 
\begin{align*}
H^i_{S_{(\pp)+}}(S_{(\pp)})_{e_0} & = 0,\ i < q,\ \text{for all $\pp \in X$},\\ H^q_{S_{(\pp)+}}(S_{(\pp)})_{e_0} & \neq  0\ \text{for some $\pp \in X$}. 
\end{align*}
It is a classical result that there exists $\dim R_{(\pp)}$ elements in $I_{(\pp)}$ which generates an ideal with the same radical as $I_{(\pp)}$. The same also holds for the ideal $S_{{(\pp)}+} = I_{(\pp)}t$. From this it follows that  
 $H_{S_{(\pp)+}}^{\dim R_{(\pp)}+1}(E) = 0$ for any $R_{(\pp)}$-module $E$ (cf. \cite[Corollary 3.3.3]{bs}.  Hence 
$$q \le \max\{\dim R_{(\pp)}|\ \pp \in X\} = \dim Y.$$
Let $Y_\pp = Y \times_X \Spec \O_{X, \pp}$. The Serre-Grothendieck correspondence yields the exact sequence
$$ 0 \to H^0_{S_{(\pp)+}}(S_{(\pp)})_{e_0} \to (S_{(\pp)})_{e_0} \to H^0(Y_\pp, \widetilde{S_{(\pp)}}(e_0)) \to H^1_{S_{(\pp)+}}(S_{(\pp)})_{e_0} \to  0, $$
and isomorphisms $H^i(Y_\pp, \widetilde{S_{(\pp)}}(e_0)) \cong H^{i+1}_{S_{(\pp)+}}(S_{(\pp)})_{e_0}$, $i \ge 1$. \par

If $q \le 1$, then $H^0(Y_\pp, \widetilde{S_{(\pp)}}(e_0)) \ne (S_{(\pp)})_{e_0} = I^{e_0}_{(\pp)}$ for some $\pp \in X$. From this it follows, as in the proof of Proposition \ref{Leray}, that
$\pi_* (\O_Y(e_0)) \neq \widetilde{I^{e_0}}.$
But $\pi_* (\O_Y(e_0))(c)$ and $\widetilde{I^{e_0}}(c)$ are generated by global sections for $c \gg 0$. Therefore, by the projection formula  we have
$$ H^0(X, \pi_* (\O_Y(c,e_0))) = H^0(X, \pi_* (\O_Y(e_0))(c)) \not= H^0(X, \widetilde{I^{e_0}}(c))\ = (I^{e_0})_c$$ 
for $c \gg 0$. Moreover,
$$H^0(Y, \O_Y(c,e_0)) = H^0(X, \pi_*(\O_Y(c,e_0))).$$
Hence $H^0(Y, \O_Y(c,e_0))  \neq (I^{e_0})_c. $\par

If $q \ge 2$, then the Serre-Grothendieck sequence implies $H^i(Y_\pp, \widetilde{S_{(\pp)}}(e_0)) = 0$  for all $\pp \in X$, $0 < i < q-1$, and 
$H^{q-1}(Y_\pp, \widetilde{S_{(\pp)}}(e_0))  \neq 0$ for some $\pp \in X$. From this it follows, as in the proof of Proposition \ref{Leray}, that
\begin{align*}
& R^i \pi_* (\O_Y(e_0)) = 0\ \text{for}\  0 < i < q-1, \\ 
& R^{q-1} \pi_* (\O_Y(e_0)) \not= 0. 
\end{align*}
By the projection formula, we have
\begin{align*}
& R^i \pi_*(\O_Y(c,e_0)) = R^i \pi_* (\O_Y(e_0)) \otimes \O_X(c) = 0\ \text{for}\ 0 < i < q-1, \\
& R^{q-1} \pi_* (\O_Y(c,e_0)) = R^{q-1} \pi_* (\O_Y(e_0)) \otimes \O_X(c) \not= 0. 
\end{align*}
Since $\pi_*(\O_Y(c,e_0)) = \pi_*(\O_Y(e_0))(c)$, we also have
$H^{q-1}(X,\pi_*(\O_Y(c,e_0))) = 0$ 
for $c \gg 0$. Therefore, using Leray spectral sequence 
$$H^i(X, R^j \pi_* (\O_Y(m,e_0))) \Rightarrow H^{i+j}(Y, \O_Y(m,e_0))$$
we can deduce that
$$H^{q-1}(Y, \O_Y(c,e_0)) = H^0(X, R^{q-1}\pi_* (\O_Y(c,e_0))).$$ 
for $c \gg 0$. But $R^{q-1} \pi_* (\O_Y(c,e_0))$ is generated by global sections for $c \gg 0$. So we get 
$H^{q-1}(Y, \O_Y(c,e_0)) \not= 0.$\par

Let us now consider the case $e_0 = a^*_X(\omega_S)$. Let $q$ be the smallest integer such that 
$e_0 = \max\{a_q((\omega_S)_{(\pp)})|\ \pp \in X\}.$ For $\pp \in X$ we have
$(\omega_S)_{(\pp)} = \oplus_{n > 0}H^0(Y_\pp,\omega_{Y_\pp}(n))$
(see \cite[2.5.2(1) and 2.6.2]{hs}). From this it follows that $[H_{S_{(\pp)+}}^i((\omega_S)_{(\pp)})]_n = 0$ for $n > 0$, $i = 0,1$. Since $e_0 > 0$, this implies $q > 1$.
Similarly as in the first case, we can also show that $q \le \dim Y$ and that $H^{q-1}(Y,\omega_Y(c,e_0)) \neq 0$ for $c \gg 0$. By Serre duality we get
$$H^{\dim Y-q+1}(Y,\O_Y(-c,-e_0)) = H^{q-1}(Y,\omega_Y(c,e_0)) \neq 0$$
for $c \gg 0$. This completes the proof of Proposition \ref{sharp-1}.
\end{proof}

We shall see later in Example \ref{sharp-example} that the bound $c > d(I)e + \max\{\varepsilon(I),\varepsilon^*(I)\}$ of Theorem \ref{diagonal-1} is sharp.\par

Now we want to study the problem when there exists a Cohen-Macaulay ring of the form $k[(I^e)_c]$ for $e \ge 1$.

\begin{Theorem} \label{diagonal-2}
Let $R$ be an equidimensional standard graded algebra over a field $k$ and $I$ a homogeneous ideal of $R$ with $\height I \ge 1$. Let $X = \Proj R$ and $S = R[It]$. Assume that $S$ is locally Cohen-Macaulay on $X$. Then, there exists a Cohen-Macaulay ring $k[(I^e)_c]$ with $c \ge d(I)e+1$ if and only if $H^0(X,\O_X) = k$ and $H^i(X,\O_X)$ $= 0$ for $i = 1,...,\dim X-1$.
Especially, this condition implies that $k[(I^e)_c]$ is a Cohen-Macaulay ring for $c > d(I)e +\max\{\varepsilon(I),\varepsilon^*(I)\}$ and $e \ge 1$.
\end{Theorem}

\begin{proof}
Let $Y$ be the blow-up of $X$ along the ideal sheaf of $I$. The assumption implies that $Y$ is equidimensional and Cohen-Macaulay. Since $S$ is locally Cohen-Macaulay over $X$, $Y$ is locally arithmetic Cohen-Macaulay over $X$. Applying Corollary \ref{transfer}, we have $H^0(Y,\O_Y) = H^0(X,\O_X)$ and $H^i(Y,\O_Y) = H^i(X,\O_X)$ for $i > 0$. Therefore, the first statement follows from Lemma \ref{criterion}. Moreover, we have $\max \{a^*_X(S), a^*_X(\omega_S)\} = 0$ by Proposition \ref{arithmetic}.
Hence the second statement follows from Theorem \ref{diagonal-1}.
\end{proof}

Note that the condition $H^0(X,\O_X) = k$ and $H^i(X,\O_X)$ $= 0$ for $i = 1,...,\dim X-1$ is satisfied if $R$ is a Cohen-Macaulay ring.\par

The following example shows that the bound $c > d(I)e +\max\{\varepsilon(I),\varepsilon^*(I)\}$ is sharp.

\begin{Example} \label{sharp-example} {\rm
Let $R = k[x_0,x_1,x_2]$ and $I = (x_1^4,x_1^3x_2,x_1x_2^3,x_2^4)$. 
It is easy to see that $S = R[It]$ is locally Cohen-Macaulay on $X = \Proj R$.
We have $I^n = (x_1,x_2)^{4n}$ for all $n \ge 2$. We have
$$a^*(I^n)  = \left\{\begin{array}{lll} 
4 & \text{if} & n = 1,\\
4n-1 & \text{if} & n \ge 2. \end{array} \right. $$
From this it follows that $\varepsilon(I) = 0$.
To compute $\varepsilon^*(I)$ we approximate $I$ by the ideal $J = (x_1,x_2)^4$.
Put $S^* = R[Jt]$. Then we have the exact sequence
$$0 \to R[It] \to R[Jt] \to k \to 0$$
From this it follows that $\omega_S = \omega_{S^*}$. 
Note that $S^*$ is a Veronese subring of the ring $T = R[(x_1,x_2)t]$ and that
$T$ is a Gorenstein ring with $\omega_T = T(-2)$. Then $\omega_{S^*} = \oplus_{n \ge 1}(x_1,x_2)^{4n-2}$. We have 
$$a(\omega_n) = a^*((x_1,x_2)^{4n-2}) = 4n-3$$
for $n \ge 1$. Hence $\varepsilon^*(I) = 0$.
By Theorem \ref{diagonal-2}, these facts imply that $k[(I^e)_c]$ is Cohen-Macaulay for $c > 4e$ and $e \ge 1$ (which can be also verified directly).
On the other hand, for $c = 4$ and $e = 1$, the ring $k[I_4] = k[x_1^4,x_1^3,x_1x_2^3,x_2^4]$ is not Cohen-Macaulay.
}
\end{Example}

There have been various criteria for the Cohen-Macaulayness of Rees algebras 
(cf. \cite{ti, huh, l, suv, aht, jk, pu}), so that one can construct various classes of ideals $I$ for which $S$ is locally Cohen-Macaulay on $X$. We list here only the most interesting applications of Theorem \ref{diagonal-2}.

\begin{Corollary} \label{intersection-1}
Let $R$ be a Cohen-Macaulay standard graded algebra over a field $k$. Let $I \subset R$ be a homogeneous ideal with $\height I \ge 1$ which is a locally complete intersection. Then $k[(I^e)_c]$ is a Cohen-Macaulay ring for all $c > d(I)e + \max\{\varepsilon(I),\varepsilon^*(I)\}$ and $e \ge 1$.
\end{Corollary}

\begin{proof} Let $X = \Proj R$.  The assumption on $I$ means that $I_\pp$ is a complete intersection ideal in $R_\pp$ for $\pp \in X$. Therefore,  $R_{(\pp)}[I_{(\pp)}t]$ is Cohen-Macaulay for all $\pp \in X$. Hence, $S = R[It]$  is locally Cohen-Macaulay on $X$. The result follows from Theorem \ref{diagonal-2}.
\end{proof}

\begin{proof} Let $X = \Proj R$.  The assumption on $I$ implies that $S = R[It]$ is locally Cohen-Macaulay on $X$.  Therefore,  the conclusion follows from Theorem \ref{diagonal-2}.
\end{proof}

\begin{Corollary} \label{smooth-1}
Let $R$ be a polynomial ring over a field $k$ of characteristic zero and $I \subset R$ a non-singular homogeneous ideal with $\height I \ge 1$. Then, $k[(I^e)_c]$ is a Cohen-Macaulay ring for  $c > d(I)e+\varepsilon(I)$ and $e \ge 1$.
\end{Corollary}

\begin{proof} The assumption implies that $I$ is locally a complete intersection. 
Hence $S = R[It]$ is locally Cohen-Macaulay on $X = \Proj R$. Let $Y = \Proj S$ Then $Y$ is a projective non-singular scheme. Let $m ,n$ be positive integers with $m \ge d(I)n+1$. Then $\O_Y(m,n)$ is a very ample invertible sheaf on $Y$ because $Y \cong \Proj k[(I^n)_m]$  \cite[Lemma 1.1]{che}. Let $\omega_S$ be the canonical module of $S$ and $\omega_Y = \widetilde{\omega_S}$. Then $H^i(Y,\omega_Y(m,n)) = 0$ for $i \ge 1$  by Kodaira's vanishing theorem. 
On the other hand,  we have
$$H^i(Y,\omega_Y(m,n)) = H^i(X,\widetilde{(\omega_S)_n}(m))$$
by Proposition \ref{Leray}. Therefore,  $H^i(X,\widetilde{(\omega_S)_n}(m)) = 0$ for $i \ge 1$. Using the Serre-Grothendieck correspondence we can deduce that $H_{R_+}^i((\omega_S)_n)_m = 0$ for $i \ge 2$.  Hence $\varepsilon^*(I) = 0$. Now, the conclusion follows from Corollary \ref{intersection-1}.
\end{proof}

\begin{Remark} When $R$ is Cohen-Macaulay, a similar result to Theorem \ref{diagonal-2} was already given by  Cutkosky and Herzog \cite[Theorem 4.1]{che}. Their result shows the existence  of a constant $\delta$ such that $k[(I^e)_c]$ is Cohen-Macaulay for $c \ge \delta e$, $e > 0$, under some assumptions on the associated graded ring $\oplus_{n \ge 0}I^n/I^{n+1}$. It is not hard to see that these assumptions imply $\max\{a_X^*(S),a_X^*(\omega_S)\} \le 0$ (see \cite[Lemma 2.1 and Lemma 2.2]{che}). Hence their result is also a consequence of Theorem \ref{diagonal-1}. Similar statements to the above two corollaries were also given in \cite{che} but without any information on the slope $\delta$.
\end{Remark}

It is not easy to compute $\varepsilon(I)$ explicitly, even when $I$ is a non-singular ideal in a polynomial ring. By a famous result of Bertram, Ein and Lazarsfeld \cite{bel} we only know that if $I$ is the ideal of a smooth complex variety cut out scheme-theoretically by hypersurfaces of degree $d_1 \ge ... \ge d_m$, then $$a_i(I^n) \le d_1n + d_2 + \cdots + d_m - \height I$$
for $ i \ge 2$ and $n \ge 1$. But we do not know any bound for $a_1(I^n)$ in terms of 
$d_1,...,d_m$. It would be of interest to find such a bound. In general, if we happen to know the minimal free resolution of $S$ over a bi-graded polynomial ring then we can estimate $\varepsilon(I)$ in terms of the shifts of syzygy modules of the resolution \cite{cht}.
\par

In the case when $I$ is the defining ideal of a scheme of fat points we know an explicit bound for $a^*(I^n)$, namely $a^*(I^n) \le \reg(I)n$ for all $n \ge 1$ \cite{ch, ggp}. As a consequence, we immediately obtain the following result of Geramita, Gimigliano and Pitteloud.

\begin{Corollary} {\rm   \cite[Theorem 2.4]{ggp})} \label{point-1}
Let $R$ be a polynomial ring over a field $k$ of characteristic zero, and $I \subset R$ the defining ideal of a scheme of fat points in $\Proj R$. Then, $k[(I^e)_c]$ is a Cohen-Macaulay ring for $c > \reg(I)e$ and $e \ge 1$.
\end{Corollary}

\begin{proof} By definition, the ideal $I$ has the form $I = \cap_{i = 1}^s \pp_i^{m_i}$, where $\pp_i$ is the defining prime ideal of a closed point in $X = \Proj R$ and $m_i \in \NN$.  Then 
$R_{(\pp)}[I_{(\pp)}t]$ is  Cohen-Macaulay for all  $\pp \in X$. In fact, we may assume that $\pp  = \pp_i$ for some $i$. Then $\pp$ is a complete intersection and  $R_{(\pp)}[I_{(\pp)}t] = R_{(\pp)}[\pp^{m_i}_{(\pp)}t]$ is a Veronese subalgebra of $R_{(\pp)}[\pp_{(\pp)}t]$. Since $R_{(\pp)}[\pp_{(\pp)}t]$ is a Cohen-Macaulay ring, so is $R_{(\pp)}[I_{(\pp)}t]$. Thus, $S = R[It]$ is locally Cohen-Macaulay on $X$. This argument also shows that $Y = \Proj S$ is smooth. Using Kodaira vanishing theorem we can show, as in the proof of Corollary \ref{smooth-1}, that $\varepsilon^*(I) = 0$. The conclusion now follows from the proof of Theorem \ref{diagonal-2} when we replace the slope $d(I)$ by $\reg(I) \ge d(I)$ and  $\varepsilon(I) $ by $0$ because of the bound $a^*(I^n) \le \reg(I)n$.
\end{proof}

It was asked in \cite{chtv} whether there exists a Cohen-Macaulay ring $k[(I^e)_c]$ for $c \gg e \gg 0$ if $R$ is a polynomial ring and $R[It]$ is Cohen-Macaulay. This question has been positively settled in \cite[Theorem 4.5]{v1}. We can make this result more precise as follows.

\begin{Corollary} 
Let $R$ be a Cohen-Macaulay standard graded algebra over a field $k$. Let $I \subset R$ be a homogeneous ideal with $\height I \ge 1$ such that $R[It]$ is Cohen-Macaulay. Then $k[(I^e)_c]$ is a Cohen-Macaulay ring for all $c > d(I)e + \max\{\varepsilon(I),\varepsilon^*(I)\}$ and $e \ge 1$.
\end{Corollary}

\section{Arithmetically Cohen-Macaulay blow-ups}

\noindent Let $X$ be a projective scheme over a field $k$. Let $\pi: Y \to X$ be the blowing up of $X$ along an ideal sheaf $\I$. We say  that $Y$ is an {\it arithmetically Cohen-Macaulay blow-up} of $X$ if there is a standard graded $k$-algebra $R$ and a homogeneous ideal $J \subset R$ with $\height J \ge 1$ such that $X = \Proj R$, $\I = \widetilde J$, and $R[Jt]$ is a Cohen-Macaulay ring. The aim of this section is to characterize arithmetically Cohen-Macaulay blow-ups. \par

Let $R$ be a finitely generated standard graded algebra over $k$, and $I$ a homogeneous ideal of $R$ with $\height I \ge 1$, such that $X = \Proj R$ and $\I = \widetilde I$. Let $d(I)$ denote the maximal degree of the elements of a homogeneous basis of $I$. For any ideal $J$ generated by $(I^e)_c$ with $c \ge d(I)e$ we have $J_n = (I^e)_n$ for all $n \ge c$ so that $\I^e = \widetilde J$. Hence $Y = \Proj R[Jt]$. The Rees algebra $R[(I^e)_ct] = R[Jt]$ is called a {\it truncated Rees algebra} of $I^e$ \cite{ha1, cha}. We may strengthen the problem on the characterization of arithmetically Cohen-Macaulay blow-ups by asking the question of when there does exist a Cohen-Macaulay truncated Rees algebra $R[(I^e)_ct]$. To solve this problem we shall need the following result of Hyry. \par

Let $T$ be a standard bi-graded algebra over a field $k$, that is, $T$ is generated over $k$ by the elements of degree $(1,0)$ and $(0,1)$. Let $M$ denote the maximal graded ideal of $T$ and define
\begin{align*}
a^1(T) & := \max\{m|\ \text{there is $n$ such that $H_M^{\dim T}(T)_{(m,n)} \neq 0$}\},\\
a^2(T) & := \max\{n|\ \text{there is $m$ such that $H_M^{\dim T}(T)_{(m,n)} \neq 0$}\}.\end{align*}\par

\begin{Theorem} \label{Hyry} \cite[Theorem 2.5]{hy}
Let $T$ be a standard bi-graded algebra over a field with $a^1(T), a^2(T) < 0$. Let $Y = \Proj T$. Then $T$ is Cohen-Macaulay if and only if the following conditions are satisfied:
\begin{align*} 
H^0(Y,T(m,n)\widetilde{~}) & \cong T_{(m,n)}\ \text{for}\ m,n \ge 0,\\
H^i(Y,T(m,n)\widetilde{~}) & = 0 \ \text{for}\ m,n \ge 0,\ i > 0,\\
H^i(Y,T(m,n)\widetilde{~}) & = 0 \ \text{for}\ m,n < 0,\ i < \dim T-2.
 \end{align*}
\end{Theorem}

Let $J \subset R$ be an arbitrary ideal generated by forms of degree $c$ and put $T = R[Jt]$. Then $T$  can be equipped with another bi-gradation given by
$$T_{(m,n)} = (J^n)_{m+cn}t^n$$
for $(m,n) \in \NN^2$. With this bi-gradation, $T$ is a standard bi-graded $k$-algebra. Comparing with the natural bi-gradation of $T$ considered in the preceding sections, we see that both bi-gradations share the same bihomogeneous elements and the same relevant bi-graded ideals. Therefore, $\Proj T$ with respect to these bi-gradations are isomorphic.\par

\begin{Lemma} \label{a*}
Let $T = R[Jt]$ be as above. Then \par
{\rm (i)} $a^1(T) \le \max\{a^*(J^n)-nc|\ n \ge 0\}$,\par
{\rm (ii)} $a^2(T) < 0$.
\end{Lemma}

\begin{proof}
To prove (i) we will show more, namely, that $H_M^i(T)_{(m,n)} = 0$ for $m  > \max\{a^*(J^n)-nc|\ n \ge 0\}$ and $i \ge 0$. Let $T_1$ denote the ideal of $T$ generated by the homogeneous elements of degree $(1,0)$. Then, by \cite[Lemma 2.3]{hy}, we only need to show that $H_{T_1}^i(T)_{(m,n)} = 0$ for $m  > \max\{a^*(J^n)-nc|\ n \ge 0\}$ and $i \ge 0$.  Since $T_1$ is generated by $R_+$, we always have
$$H_{T_1}^i(T)_{(m,n)} = \left\{\begin{array}{lll}
0 & \text{for} & n < 0,\\
H_{R_+}^i(J^n)_{m+nc} & \text{for} & n \ge 0.
\end{array}\right.$$
But $H_{R_+}^i(J^n)_{m+nc} = 0$  for $m + nc > a^*(J^n)$, $n \ge 0$. Therefore, $H_{T_1}^i(T)_{(m,n)} = 0$ for $m  > \max\{a^*(J^n)-nc|\ n \ge 0\}$, as required.\par
To prove (ii) we first observe that 
$$a^2(T) = \max\{n|\ H_M^{\dim T}(T)_n \neq 0\},$$
where the $\ZZ$-gradation comes from the natural grading 
$T_n = J^nt^n$, $n \ge 0$. Therefore, the conclusion $a^2(T) < 0$ follows from \cite[Corollary 3.2]{ti}.
\end{proof}

\begin{Corollary} \label{a(T)}
Let $R$ be a standard graded algebra over a field with $a^*(R) < 0$ and $I \subset R$ a homogeneous ideal with $\height I \ge 1$. Let  $T = R[(I^e)_ct]$ for some fixed integers $c > d(I)e+ \varepsilon(I)$ and $e \ge 1$. Then
$a^1(T) < 0$ and $a^2(T) < 0$.
\end{Corollary}

\begin{proof}
Let $J$ be the ideal of $R$ generated by $(I^e)_c$. By Lemma \ref{a*} we only need to prove that $a^*(J^n) < nc$ for $n \ge 0$. For $n = 0$, this follows from the assumption $a^*(R) < 0$. For $n \ge 1$, we will approximate $a^*(J)$ by $a^*(I^{en})$. Since $J^n$ is generated by elements of degree $cn$ and since 
$cn > d(I)en \ge d(I^{en})$, we have $(I^{en}/J^n)_m = 0$ for $m \ge cn$. From this it follows that $H^0(I^{en}/J^n) = I^{en}/J^n$ and $H^i(I^{en}/J^n) = 0$ for $i > 0$. Therefore, from the  exact sequence 
$$0 \To J^n \To I^{en} \To I^{en}/J^n \To 0$$
we can deduce that $H^i(J^n)_m = H^i(I^{en})_m$ for $m \ge cn$ and $i \ge 0$. This implies 
$$a^*(J^n) \le \max\{cn-1,a^*(I^{en})\}.$$
By the definition of $\varepsilon(I)$ we have $a^*(I^{en}) \le d(I)en + \varepsilon(I) \le cn -1.$ Therefore, $a^*(J^n) \le cn-1$ for $n \ge 1$. \end{proof}

We are now ready to give a necessary and sufficient condition for the existence of a Cohen-Macaulay truncated Rees algebra. 

\begin{Theorem} \label{truncated-1}
Let $R$ be a standard graded algebra over a field with $a^*(R) < 0$ and $I \subset R$ a homogeneous ideal with $\height I \ge 1$. Let $X = \Proj R$, $S = R[It]$ and $Y = \Proj S$. Then there exists a Cohen-Macaulay ring $R[(I^e)_ct]$ with $c \ge d(I)e$ if and only if  the following conditions are satisfied:\par
{\rm (i)} $Y$ is equidimesional and Cohen-Macaulay,\par
{\rm (ii)} $\pi_*\O_Y = \O_X$ and $R^i\pi_*\O_Y = 0$ for $i > 0$.\par
\noindent Especially, these conditions imply that $R[(I^e)_ct]$ is a Cohen-Macaulay ring for $c > d(I)e+\max\{\varepsilon(I),\varepsilon^*(I)\}$ and $e > \max\{a_X^*(S),a_X^*(\omega_S)\}$. 
\end{Theorem}

\begin{proof}
Let $J$ be the ideal of $R$ generated by $(I^e)_c$ and $T = R[Jt]$ for a fixed pair of positive integers $c,e$ with $c \ge d(I)e$. Then $Y \cong \Proj T$. 
If $T$ is a Cohen-Macaulay ring, then (i) is obviously satisfied and $Y$ is locally arithmetic Cohen-Macaulay over $X$. (ii) follows from Corollary \ref{transfer}. \par

To prove the converse we equip $T$ with the afore mentioned bi-gradation.
Set $e_0 = \max\{a_X^*(S),a_X^*(\omega_S)\}$. 
We will use Theorem \ref{Hyry} to prove that $T$ is Cohen-Macaulay for $c > d(I)e+\max\{\varepsilon(I),\varepsilon^*(I)\}$ and $e > e_0$. 
By Corollary \ref{a(T)} we have $a^1(T) < 0$ and $a^2(T) < 0$. From the bi-gradation of $T$ we see that 
$$T(m,n)\widetilde{~} = \O_Y(m+cn,en),$$
where $\O_Y(m+cn,n)$ denotes the twisted $\O_Y$-module with respect to the natural bi-gradation of $S$. If $\pi_*\O_Y = \O_X$ and $R^i\pi_*\O_Y = 0$ for $i > 0$, then we can show as in the proof of Proposition \ref{Leray} that
$H^i(Y,\O_Y(m,0))
 = H^i(X,\O_X(m))$
for $i \ge 0$. Since $a^*(R) < 0$, we have $H_{R_+}^i(R)_m = 0$ for all $m \ge 0$ and $i \ge 0$.
Using the Serre-Grothendieck correspondence between sheaf cohomology of $X$ and local cohomology of $R$ we can deduce that $H^0(X,\O_X(m)) = R_m$ and $H^i(X,\O_X(m)) = 0$ for $i > 0$. Therefore,
\begin{align*}
H^0(Y, \O_Y(m,0)) & = R_m = T_{(m,0)},\\
H^i(Y, \O_Y(m,0)) & = 0, \ i > 0. 
\end{align*}
For $m \ge 0$ and $n > 0$ we have $m +cn > d(I)en + \varepsilon(I)$. Therefore, using Proposition \ref{Leray} and Lemma \ref{asymptotic} we get 
\begin{align*}
H^0(Y,\O_Y(m+cn,en))  & =  T_{(m,n)},\\
H^i(Y,\O_Y(m+cn,en))  & =  0,\   i > 0,
\end{align*}
for $e > e_0$. For $m, n < 0$ we can show, similarly as above, that
$H^i(Y,\omega_Y(-m-cn,-en)) = 0$ for $i > 0$ and $e > e_0$. If $Y$ is equidimensional and Cohen-Macaulay, we can apply Serre duality and obtain 
$$H^i(Y,\O_Y(m+cn,en)) = 0,\  i < \dim Y.$$
Passing from $\O_Y(m+cn,en)$ to $T(m,n)\widetilde{~}$ we get
\begin{align*} 
H^0(Y,T(m,n)\widetilde{~}) & \cong T_{(m,n)}\ \text{for}\ m,n \ge 0,\\
H^i(Y,T(m,n)\widetilde{~}) & = 0 \ \text{for}\ m,n \ge 0,\ i > 0,\\
H^i(Y,T(m,n)\widetilde{~}) & = 0 \ \text{for}\ m,n < 0,\ i < \dim T-2.
\end{align*}
By Theorem \ref{Hyry}, these conditions imply that $T$ is a Cohen-Macaulay ring.
The proof of Theorem \ref{truncated-1} is now complete.
\end{proof}

The following example shows that the condition $a^*(R) < 0$ is not necessary for the existence of a Cohen-Macaulay truncated Rees algebra. It also shows that in general, the existence of a  Cohen-Macaulay truncated Rees algebra does not imply the existence of a linear bound on $c$ ensuring the Cohen-Macaulayness of $R[(I^e)_ct]$.

\begin{Example} \label{non-linear} Take $R = k[x,y,z]/(xy^2-z^3)$, the coordinate ring of a plane cusp, and $I = (x) \subseteq R$, a homogeneous ideal with $\height I = 1$. Then 
$R$ is a two-dimensional Cohen-Macaulay ring with  $a^*(R) = 0$. It is obvious that $R[(I^e)_et] = R[It]$ is a Cohen-Macaulay ring for $e \ge 1$.
For $c > e$ we have 
$R[(I^e)_ct] \cong R[(x,y,z)^{c-e}t]$. It is easy to check that the 
reduction number of the ideal $(x,y,z)^{c-e}$ is greater than 1. By \cite{gs}, this implies that $R[(x,y,z)^{c-e}t]$ is not Cohen-Macaulay for any $c > e$.
\end{Example}

Now we will show that the bound $e > e_0$ in Theorem \ref{truncated-1} is once again best possible.

\begin{Proposition} \label{sharp-2}
Let the notations and assumptions be as in Theorem \ref{truncated-1}. Put
\[ e_0 = \max \{ a^*_X(S), a^*_X(\omega_S) \}. \]
Then $R[(I^{e_0})_c t]$ is not a Cohen-Macaulay ring for  $c \ge d(I)e_0$ if $e_0 \ge 1$.
\end{Proposition}

\begin{proof} 
Let $e_0 \ge 1$ and $T = R[(I^{e_0})_c t]$ for some $c \ge d(I)e_0$. Note that $(I^{e_0})_c$ and $I^{e_0}$ defines the same ideal sheaf in $\O_X$. Consider the natural $\NN$-grading of $T$ and $S$ given by the degree of $t$. For any $\pp \in X$, the ring $T_{(\pp)}$ is isomorphic to the $e_0$-th Veronese subring of $S_{(\pp)}$.  Hence
\begin{align*}
& H_{T_{(\pp)+}}^i(T_{(\pp)})_1 = H_{S_{(\pp)+}}^i(S_{(\pp)})_{e_0},\\
& H_{T_{(\pp)+}}^i((\omega_T)_{(\pp)})_1 = H_{S_{(\pp)+}}^i((\omega_S)_{(\pp)})_{e_0},
\end{align*}
for $i \ge 0$. By the definition of $e_0$ there exists $\pp \in X$ and $i \ge 0$ such that either $H_{S_{(\pp)+}}^i(S_{(\pp)})_{e_0} \neq 0$ or $H_{S_{(\pp)+}}^i((\omega_S)_{(\pp)})_{e_0} \neq 0$. Therefore, $\max\{a^*(T),a^*(\omega_T)\} \ge 1$. By Corollary \ref{transfer}, this implies that $T$ is not a Cohen-Macaulay ring.
\end{proof}

From Theorem \ref{truncated-1} we can derive the following sufficient  condition for the existence of a truncated Cohen-Macaulay Rees algebra.

\begin{Theorem} \label{truncated-3}
Let $R$ be an equidimensional standard graded algebra over a field with $a^*(R) < 0$ and $I \subset R$ a homogeneous ideal with $\height I \ge 1$. Let $X = \Proj R$ and $S = R[It]$. Assume that $S$ is locally Cohen-Macaulay on $X$. Then $R[(I^e)_ct]$ is a Cohen-Macaulay ring for $c > d(I)e+\max\{\varepsilon(I),\varepsilon^*(I)\}$ and $e \ge 1$. 
\end{Theorem}

\begin{proof}
It is obvious that the assumptions imply that $Y$ is equidimensional and Cohen-Macaulay. The condition $\pi_*\O_Y = \O_X$ and $R^i\pi_*\O_Y = 0$ for $i > 0$ follows from Corollary \ref{transfer}. Hence the conclusion follows from Theorem \ref{truncated-1}.
\end{proof}

The above condition is also a necessary condition for the existence of a truncated Cohen-Macaulay Rees algebra of the form $R[I_ct]$ ($e = 1$).

\begin{Corollary} \label{e=1}
Let $R$ be a standard graded algebra over a field with $a^*(R) < 0$ and $I \subset R$ a homogeneous ideal with $\height I \ge 1$. Let $X = \Proj R$ and $S = R[It]$. Then there exists a Cohen-Macaulay ring $R[I_ct]$ with $c \ge d(I)$ if and only if  $S$ is locally Cohen-Macaulay on $X$. 
\end{Corollary}

\begin{proof}
By Theorem \ref{truncated-3} we only need to show that if $R[I_ct]$ is a Cohen-Macaulay ring for some $c \ge d(I)$, then $S$ is locally Cohen-Macaulay on $X$.
But this is obvious because $(I_c)$ and $I$ define the same ideal sheaf and $R[I_ct]$ is locally Cohen-Macaulay on $X$. 
 \end{proof}

Using Theorem \ref{truncated-3} we obtain several classes of Cohen-Macaulay Rees algebras.

\begin{Corollary} \label{intersection-2}{\rm (cf. \cite[Corollary 2.2.1(2)]{cha} for the case $e=1$)}
Let $R$ be a Cohen-Macaulay standard graded algebra over a field $k$ with $a(R) < 0$. Let $I \subset R$ be a homogeneous ideal with $\height I \ge 1$ which is locally a complete intersection. Then  $R[(I^e)_ct]$ is a Cohen-Macaulay ring for all $c > d(I)e + \max\{\varepsilon(I),\varepsilon^*(I)\}$ and $e \ge 1$.
\end{Corollary}

\begin{proof} As in the proof of Corollary \ref{intersection-1}, $S = R[It]$ is locally Cohen-Macaulay over $X = \Proj R$.  Since the assumption on $R$ implies $a^*(R) < 0$, the conclusion follows from Theorem \ref{truncated-3}.
\end{proof}

\begin{Corollary} \label{smooth-2}
Let $R$ be a polynomial ring over a field $k$ of characteristic zero  and $I \subset R$ a non-singular homogeneous ideal. Then  $R[(I^e)_ct]$ is a Cohen-Macaulay ring for all $c > d(I)e+\varepsilon(I)$ and $e \ge 1$.
\end{Corollary}

\begin{proof}
We have seen in the proof of Corollary \ref{smooth-1} that $\varepsilon^*(I) = 0$. Hence the assertion follows from Corollary \ref{intersection-2}.
\end{proof}

\begin{Corollary} {\rm (cf. \cite[Theorem 2.4]{ha1} for the case $e=1$)}
Let $R$ be a polynomial ring over a field $k$ of characteristic zero and $I \subset R$ the defining ideal of a scheme of fat points in $\Proj R$. Then $R[(I^e)_ct]$ is a Cohen-Macaulay ring for $c > \reg(I)e$.
\end{Corollary}

\begin{proof} The proof follows from Theorem \ref{truncated-3} with the same lines of arguments as in the proof of Corollary \ref{point-1}.
\end{proof}

Now we will use Theorem \ref{truncated-3} to find a criterion for  arithmetically Cohen-Macaulay blow-ups. Recall that the blow-up $Y$ of a projective scheme $X$ along an ideal sheaf $\I$ is said to be {\it locally arithmetic Cohen-Macaulay on $X$} if there exist a standard graded algebra $R$ over a field and a homogeneous ideal $I \subset R$ such that $X = \Proj R$, $\I = \widetilde{I}$ and $S = R[It]$ is locally Cohen-Macaulay on $X$. 

\begin{Theorem} \label{truncated-2}
Let $X$ be a projective scheme over a field $k$ such that $H^0(X,\O_X) = k$ and $H^i(X,\O_X) = 0$ for $i > 0$. Let $Y$ be a blow-up of $X$. Then $Y$ is an arithmetically Cohen-Macaulay blow-up if and only if $Y$ is equidimensional and locally arithmetic Cohen-Macaulay on $X$.
\end{Theorem}

\begin{proof} Suppose $Y$ is an arithmetically Cohen-Macaulay blow-up of $X$. Let $R$ be a standard graded algebra over $k$, and $I$ be a homogeneous ideal of $R$, such that $X = \Proj R$, $Y$ is the blow-up of $X$ along the ideal sheaf $\widetilde I$, and $S = R[It]$ is a Cohen-Macaulay ring. Then, $\O_{X,x}[\I_xt] = S_{(\pp)}$ is obviously Cohen-Macaulay for all $\pp \in X$. Thus, $Y$ is locally arithmetic Cohen-Macaulay on $X$. 

Conversely, suppose $Y$ is equidimensional and locally arithmetic Cohen-Macaulay on $Y$. Then there exist a standard graded $k$-algebra $R$ and a homogeneous ideal $I \subset R$ such that $X = \Proj R$, $Y$ is the blow-up of $X$ along the ideal sheaf of $I$, and $R[It]$ is locally Cohen-Macaulay on $X$. The assumption on the sheaf cohomology of $X$ implies that $H_{R_+}^{i}(R)_0 = 0$ for $i \ge 0$.  Without restriction we may replace $R$ by a suitable Veronese subalgebra and obtain $H_{R_+}^{i}(R)_n = 0$ for all $n \ge 0$ or, equivalently, $a^*(R) < 0$. Now we may apply Theorem \ref{truncated-3} to find a Cohen-Macaulay Rees algebra $R[I_ct]$ with $c \gg 0$. Since the ideal $(I_c)$ defines the same ideal sheaf $\widetilde I$, we can conclude that $Y$ is an arithmetically blow-up of $X$.
\end{proof}

%%%%%%%%%%%%%%%%%%%%%%%%%%%%%%%%%%%%%%%%%%%%%%%%%%%%%%%%%%%%%%%%%%%


\begin{thebibliography}{100}
\bibitem{aht} I. M. Aberbach, C. Huneke and N.V. Trung. {\it Reduction numbers, Briancon-Skoda theorems, and depth of Rees rings}. Compositio Math. {\bf 97} (1995), 403-434.
\bibitem{ao} Y. Aoyama. {\it On the depth and the projective dimension of the canonical module}. Japanese J. Math. {\bf 6} (1980), 61-66.
\bibitem{bel} A. Bertram, L. Ein and R. Lazarsfeld. {\it Vanishing
theorems, a theorem of Severi, and the equations defining
projective varieties}. J. Amer. Math. Soc. {\bf 4} (1991), no. 3, 587-602.
\bibitem{bs} M. Brodmann and R. Sharp.  Local cohomology. Cambridge University Press, 1998. 
\bibitem{bh} W. Bruns and J. Herzog. Cohen-Macaulay rings. Cambridge University Press, 1993.
\bibitem{ch} K. A. Chandler. {\it Regularity of the powers of an
ideal}. Commun. Algebra. {\bf 25} (1997), 3773-3776.
\bibitem{chtv} A. Conca, J. Herzog, N.V. Trung and G. Valla. {\it Diagonal subalgebras of bi-graded algebras and embeddings of blow-ups of projective spaces}. American Journal of Math. {\bf 119} (1997), 859-901.
\bibitem{cha} S.D. Cutkosky and H. T\`ai H\`a. {\it Arithmetic Macaulayfication of projective schemes}. J. Pure Appl. Algebra. To appear. \bibitem{che} S.D. Cutkosky and J. Herzog. {\it Cohen-Macaulay coordinate rings of blowup schemes}. Comment. Math. Helv. {\bf 72} (1997), 605-617.
\bibitem{cht} S.D. Cutkosky, J. Herzog and N.V. Trung. {\it Asymptotic
behaviour of the Castelnuovo-Mumford regularity}. Compositio Math. {\bf 118} (1999), 243-261.
\bibitem{gg} A.V. Geramita and A. Gimigliano. {\it Generators for the defining ideal of certain rational surfaces}. Duke Mathematical Journal. {\bf 62} (1991), no. 1, 61-83.
\bibitem{ggh} A.V. Geramita, A. Gimigliano and B. Harbourne. {\it Projectively normal but superabundant embeddings of rational surfaces in projective space}. J. Algebra. {\bf 169} (1994), no. 3, 791-804.
\bibitem{ggp} A.V. Geramita, A. Gimigliano and Y. Pitteloud. {\it Graded Betti numbers of some embedded rational $n$-folds}. Math. Ann. {\bf 301} (1995), 363-380.
\bibitem{gs} S. Goto and Y. Shimoda, {\it On the Rees algebras of Cohen-Macaulay rings}. Lect. Notes in Pure and Appl. Math. 68, Marcel-Dekker, 1979, 201-231.
\bibitem{ha1} H. T\`ai H\`a. {\it On the Rees algebra of certain codimension two perfect ideals}. Manu. Math. {\bf 107} (2002), 479-501.
\bibitem{ha2} H. T\`ai H\`a. {\it Projective embeddings of projective schemes blown up at subschemes}. Math. Z. To appear.
\bibitem{har} R. Hartshorne. Algebraic Geometry. Graduate Text {\bf 52}. Springer-Verlag, 1977.
\bibitem{huh} S. Huckaba and C. Huneke. {\it Rees algebras of ideals having small analytic deviation}. Trans. Amer. Math. Soc. {\bf 339} (1993), no. 1, 373-402.
\bibitem{hy} E. Hyry. {\it The diagonal subring and the Cohen-Macaulay property of a multigraded ring}. Trans. Amer. Math. Soc. {\bf 351} (1999), no. 6, 2213-2232.
\bibitem{hs} E. Hyry and K. Smith. {\it On a Non-Vanishing Conjecture of Kawamata and the Core of an Ideal}. Preprint. {\tt arXiv:math.AG/0301189}
\bibitem{jk} B. Johnston and D. Katz. {\it Castelnuovo regularity and graded rings associated to an ideal}. Proc. Amer. Math. Soc. {\bf 123} (1995), 727-734.
\bibitem{ka} T. Kawasaki. {\it On arithmetic Macaulayfication of local rings}. Trans. Amer. Math. Soc. {\bf 354}, 123-149.
\bibitem{k} V. Kodiyalam. {\it Asymptotic behaviour of
Castelnuovo-Mumford regularity}. Proc. Amer. Math. Soc. {\bf 128}
(2000), 407-411.
\bibitem{m} D. Mumford. {\it Varieties defined by quadratic equations}. C.I.M.E. III. (1969), 29-100.
\bibitem{v1} O. Lavila-Vidal. {\it On the Cohen-Macaulay property of diagonal subalgebras of the Rees algebra}. Manu. Math. {\bf 95} (1998), 47-58.
\bibitem{v2} O. Lavila-Vidal. {\it On the existence of Cohen-Macaulay coordinate rings of blow-up schemes}. Preprint.
\bibitem{l} J. Lipman. {\it Cohen-Macaulayness in graded algebras}. Math. Res. Letters {\bf 1} (1994), 149-157.
\bibitem{pu} C. Polini and B. Ulrich. {\it Neccessary and sufficient conditions for the Cohen-Macaulayness of blow-up algebras}. Compositio Math. {\bf 119} (1999), no. 2, 185-207.
\bibitem{sh} R. Sharp. {\it Bass numbers in the graded case, $a$-invariant formula, and an analogue of Falting's annihilator theorem}. J. Algebra. {\bf 222} (1999), no. 1, 246-270.
\bibitem{suv} A. Simis, B. Ulrich, and W. Vasconcelos. {\it Cohen-Macaulay Rees algebras and degrees of plolynomial equations}. Math. Ann. {\bf 301} (1995), 421-444.
\bibitem{sw} I. Swanson. {\it Powers of ideals. Primary
decompositions, Artin-Rees lemma and regularity}. Math. Ann. {\bf 307}
(1997), 299-313.
\bibitem{t1} N.V. Trung. {\it The largest non-vanishing degree of graded local cohomology modules}. J. Algebra. {\bf 215} (1999), no. 2, 481-499.
\bibitem{ti} N.V. Trung and S. Ikeda. {\it When is the Rees algebra Cohen-Macaulay?} Comm. Algebra. {\bf 17} (1989), no. 12, 2893-2922. 
\bibitem{tw} N.V. Trung and H-J. Wang. {\it On the asymptotic linearity of Castelnuovo-Mumford regularity}. Preprint. {\tt arXiv:math.AC/0212161}.
\end{thebibliography}
\end{document}